\documentclass[a4paper,twoside]{amsart}
\usepackage{amsmath,amsfonts,amssymb,amscd,amsthm,latexsym,enumerate}
\usepackage[all]{xy}
\newcounter{weiter}
\def\polhk#1{\setbox0=\hbox{#1}{\ooalign{\hidewidth
    \lower1.5ex\hbox{`}\hidewidth\crcr\unhbox0}}}

\newtheorem{Lem}{Lemma}[section]
\newtheorem{Prop}[Lem]{Proposition}
\newtheorem{Cor}[Lem]{Corollary}
\newtheorem{Thm}[Lem]{Theorem}

\theoremstyle{definition}
\newtheorem{Def}[Lem]{Definition}

\newtheorem{Rem}[Lem]{Remark}

\renewcommand\o{\otimes}

\newcommand\op{{\operatorname{op}}}
\newcommand\cop{{\operatorname{cop}}}

\newcommand\entwmod[3]{{\mathcal M_{#1}^{#2}}}
\newcommand\HMod[4]{{^{#1}_{#3}\mathcal M^{#2}_{#4}}}

\def\HM#1.#2.#3.#4.{\HMod{#1}{#2}{#3}{#4}}

\newcommand\LMod[1]{{_{#1}\mathcal M}}
\newcommand\LComod[1]{{^{#1}\mathcal M}}
\newcommand\RComod[1]{\mathcal M^{#1}}

\newcommand\RMod[1]{\mathcal M_{#1}}

\newcommand\lcofix[2]{{^{\operatorname{co} #2}{#1}}}
\newcommand\rcofix[2]{{#1}{^{\operatorname{co} #2}}}
\newcommand{\leer}{\operatorname{--}}
\newcommand{\ou}[1]{\mathrel{\mathop{\otimes}_{#1}}}
\newcommand{\co}[1]{\mathrel{\mathop{\Box}_{#1}}}

\newcommand\hitby{\leftharpoonup}

\newcommand\sw[1]{{}_{(#1)}}
\newcommand\swm[1]{{}_{(-#1)}}

\newcommand\ol{\overline}
\newcommand\inv{^{-1}}
\renewcommand\epsilon\varepsilon

\makeatletter
\def\namelabel#1#2{\@bsphack
  \protected@write\@auxout{}%
         {\string\newlabel{#1.nme}{{#2}{#2}}}%
  \@esphack}
\def\nmlabel#1#2{\label{#2}\namelabel{#2}{#1}}
\newcommand\nmref[1]{\ref{#1.nme}\ \ref{#1}}
\makeatother

\newcommand\bpf{\begin{proof}}
\newcommand\epf{\end{proof}}
\def\Label#1{\label{#1}\ifmmode\llap{[#1] }\else
\marginpar{\smash{\hbox{[#1]}}}\fi}

\begin{document}


\newtheorem{proposition}{\sc Proposition}[section]
\newtheorem{lemma}[proposition]{\sc Lemma}
\newtheorem{corollary}[proposition]{\sc Corollary}
\newtheorem{theorem}[proposition]{\sc Theorem}

\theoremstyle{definition}
\newtheorem{definition}[proposition]{\sc Definition}
\newtheorem{example}[proposition]{\sc Example}

\theoremstyle{remark}
\newtheorem{remark}[proposition]{\sc Remark}

\def\proof{{\sc Proof.~~}}
\def\endproof{\hbox{$\sqcup$}\llap{\hbox{$\sqcap$}}\vskip 4pt plus2pt}

\newcommand{\Section}{\setcounter{definition}{0}\section}
\renewcommand{\theequation}{\thesection.\arabic{equation}}
\newcounter{c}
\renewcommand{\[}{\setcounter{c}{1}$$}
\newcommand{\etyk}[1]{\vspace{-7.4mm}$$\begin{equation}\Label{#1}
\addtocounter{c}{1}}
\renewcommand{\]}{\ifnum \value{c}=1 $$\else \end{equation}\fi}
\setcounter{tocdepth}{2}


\newcommand{\can}{\operatorname{can}}
\newcommand{\alg}{\operatorname{Alg}}
\newcommand{\id}{\operatorname{id}}
\newcommand{\Me}{\mathcal{M}_A^C(\psi)}
\newcommand{\Mp}{\mathcal{M}_P^C(\psi)}
\newcommand{\Es}{(A,C)_{\psi}}
\newcommand{\Ep}{(P,C)_{\psi}}
\newcommand{\Op}{\Omega^1 P}
\newcommand{\Ob}{\Omega^1 B}
\newcommand{\M}{\mathcal{M}}
\newcommand{\C}{\mathcal{C}}
\newcommand{\Ker}{\operatorname{Ker}}
\newcommand{\Ima}{\operatorname{Im}}
\newcommand{\Hom}{\operatorname{Hom}}
\newcommand{\tr}{\operatorname{Tr}}
\def\sw#1{{\sb{(#1)}}}
\def\sco#1{{\sp{(#1)}}} 
\def\su#1{{\sp{[#1]}}} 
\def\eps{\varepsilon}
\def\DC{\Delta_{\C}}
\def\eC{\eps_\C}
\def\ut{\otimes}
\def\ov{\overline}
\def\act{\triangleleft}
\def\csop{\mbox{$(C^*)^{op}$}}
\def\lfa{{\large\mbox{$\forall\;$}}}
\def\fa{\forall\,\,}

\renewcommand{\labelenumi}{\theenumi)}
\newcommand{\nc}[2]{\newcommand{#1}{#2}}
\newcommand{\rnc}[2]{\renewcommand{#1}{#2}}

\rnc{\theequation}{\thesection.\arabic{equation}}

\hyphenation{Gro-then-dieck}

\newtheorem{examples}[proposition]{\sc Examples}

\nc{\bexs}{\begin{examples}}
\nc{\eexs}{\hfill\mbox{$\Diamond$}\end{examples}}

\newenvironment{enumerate1}{\begin{enumerate}}
{\end{enumerate}}
\renewcommand\labelenumi{\rm(\arabic{enumi})}
\renewcommand\labelenumii{\rm(\alph{enumii})}
\title{On generalized Hopf Galois extensions}
\author{Peter Schauenburg
and Hans-J\"urgen Schneider}
\thanks{Peter Schauenburg thanks the DFG for support by a Heisenberg Fellowship}
\address{Mathematisches Institut der Universit\"at
M\"unchen\\Theresienstr.~39\\80333 M\"unchen\\Germany}
\email{hanssch@mathematik.uni-muenchen.de,schauen@mathematik.uni-muenchen.de}
\maketitle
\section*{Introduction}
Let $H$ be a Hopf algebra over a commutative base ring $k$, and
$A$ a right $H$-comodule algebra with comodule structure
$\delta\colon A\to A\o H$, $\delta(a)=:a\sw 0\o a\sw 1$. Denote by
$B:=\rcofix AH:=\{b\in A|\delta(b)=b\o 1\}$ the subalgebra of
coinvariant elements. $A$ is said to be an $H$-Galois extension of
$B$ if the Galois map
$$\beta\colon A\ou BA\to A\o H,x\o y\mapsto xy\sw 0\o y\sw 1,$$
is a bijection.

A faithfully flat (as $B$-module) $H$-Galois extension $A$ is a
noncommutative-geometric version of a principal fiber bundle or
torsor in the sense of \cite{DemGab:GAI}: If $A$ and $H$ are
commutative, and represent respectively an affine scheme $X$ and
an affine group scheme $G$ acting on $X$, then $B=\rcofix AH$
represents the quotient $Y$ of $X$ under the action of $G$.
Bijectivity of the Galois map $\beta$ means that
$$X\times G\to X\times_YX,(x,g)\mapsto (xg,x),$$
is an isomorphism, which can be interpreted as the correct
algebraic formulation of the condition that the $G$-action of $X$
should be free, and transitive on the fibers of the map $X\to Y$.

In many applications surjectivity of the Galois map
$\beta$, which, in the
commutative case, means freeness of the action of $G$, 
is obvious, or at least easy to prove (it is sufficient to
find $1\o h$ in the image for each $h$ in a generating set for the
algebra $H$). It is usually much harder to decide whether $\beta$ is
injective.

The present paper has two main topics: When does surjectivity of
$\beta$ already imply bijectivity? What can we conclude about the
module structure of $A$ over $B$, or the comodule structure of 
$A$, or general Hopf modules, over $H$? Both questions will be 
studied for more general extensions.

The Kreimer-Takeuchi Theorem \cite[Thm.~1.7]{KreTak:HAGEA} says
that if $\beta$ is onto and $H$ is finite, then $\beta$ is
bijective and $A$ is a projective $B$-module. This generalizes a
Theorem of Grothendieck \cite[III, \S 2, 6.1]{DemGab:GAI} on the
actions of finite group schemes. Theorem 3.5 in \cite{Sch:PHSAHA}
implies that if  $\beta$ is surjective and $A$ is a relative
injective $H$-comodule, then $\beta$ is bijective and $A$ is a
faithfully flat $B$-module. This generalizes results of Oberst
\cite{Obe:AQAAGID}, and Cline, Parshall, and Scott
\cite{CliParSco:IMAQ} for the case where $H$ represents a closed
subgroup of an affine group scheme represented by $A$; in this
situation the canonical map is trivially surjective, while
injectivity of the $H$-comodule $A$ means that the induction
functor from the subgroup in question is exact.

A new proof for both of these results appeared in
\cite{Sch:HGBGE}, where it is also shown that in the situation of
\cite[Thm.I]{Sch:PHSAHA} the $B$-module $A$ is projective as well.
The unified proof and the stronger conclusion are based on the
observation that the Galois map $\beta_0\colon A\o A\to A\o H$
(where the tensor product in the source is taken over $k$ rather
than the coinvariant subalgebra $B$), which is surjective by
assumption, can be shown to be split as an $H$-comodule map in
each case.

In the present paper we will show (with a further simplified
proof) that having an $H$-colinearly split surjective map $\beta_0$
characterizes relative projective $H$-Galois extensions. We will
in fact show this for more general extensions, and we will discuss
applications of the generalized result, with appropriate
additional hypotheses, to a variety of Galois-type situations.

In \nmref{sec:inj} we give a new proof for
the criteria \cite[Thm.3.5, Thm.I]{Sch:PHSAHA} mentioned above.

In \nmref{sec:KTTT} we prove a strong generalization of the
Kreimer-Takeuchi Theorem. Among its corollaries are the original
Kreimer-Takeuchi Theorem as well as a result of 
Beattie, D{\u{a}}sc{\u{a}}lescu, and Raianu
\cite{BeaDasRai:GECFHA}
for the co-Frobenius case; we improve on the latter by proving 
that the extension is projective rather than flat.

In \nmref{strongconnectiontype} we will consider
another condition on a Hopf Galois extension, which we call
equivariant projectivity. This is a stronger requirement than
projectivity of $A$ over $B$; it was studied by
D{\polhk{a}}browski, Grosse, and Hajac \cite{DabGroHaj:SCCCPHGT},
who showed that a Hopf Galois extension is equivariantly
projective if and only if it has a so-called strong connection.
This notion in turn was defined by Hajac \cite{Haj:SCQPB} with
motivations from differential geometry; see also
\cite{HajMaj:PMDM}. Most notably we will show in
\nmref{PropHopfGal} that if $H$ is a Hopf algebra with bijective
antipode over a field, then every faithfully flat $H$-Galois
extension is equivariantly projective. Thus, strong connections
always exist in the situation for which they were originally
defined.

The reason why we are interested in generalizations of Hopf Galois
extensions lies in the quotient theory of noncommutative Hopf
algebras. The quotient Hopf algebras of a commutative Hopf algebra
$H$ correspond naturally to the closed subgroups of the affine
group scheme represented by $H$. If $H$ is noncommutative,
however, it is not enough to consider quotient Hopf algebras.
Rather, one should also take into account quotient coalgebras and
right (or left) $H$-modules $Q$ of $H$; quotient theory of Hopf 
algebras in this sense was studied by Takeuchi \cite{Tak:QSHA} and 
Masuoka \cite{Mas:QTHA}. Thus it becomes natural to
consider $Q$-Galois extensions, that is, $H$-comodule algebras $A$
for which the canonical map $A\ou BA\to A\o Q$ is bijective, for
$B=\rcofix AQ$. Such extensions are already studied in
\cite{Sch:PHSAHA}. It is important for the theory that the notion
of a Hopf module $M\in\HM.Q..A.$ can be defined for a right
$H$-module coalgebra quotient $Q$ of $H$ in the same way as a Hopf
module in $\HM.H..A.$. Later $Q$-extensions were studied in
successively more general frameworks. In the most general version
a Coalgebra Galois extension \cite{BrzHaj:CEACGT} is simply an
algebra $A$ which is a $C$-comodule for a coalgebra $C$ such that
the canonical map $\beta\colon A\ou BA\to A\o C$ is a bijection;
here $B=\rcofix AC$ is defined, following Takeuchi, as the largest
subalgebra for which the comodule structure of $A$ is left
$B$-linear. The notion of a Hopf module, which is a central tool
for studying Hopf Galois extensions, also underwent a series of
generalizations: First, one can study an $H$-comodule algebra $A$
as before, but replace $Q$ by an $H$-module coalgebra $C$; thus,
one arrives at the self-dual notion of Doi-Koppinen data $(A,C,H)$
for which much of the theory of Hopf modules can be developed
\cite{Doi:UHM,Kop:VSPAGGR}. The Hopf algebra's main role in this
formalism is to induce a generalized switching map
$$C\o A\ni c\o a\mapsto a\sw 0\o c\cdot a\sw 1\in A\o C$$
between the algebra and coalgebra in consideration. A further step
abstracts this switching map from the auxiliary Hopf algebra $H$
and defines an entwining \cite{BrzMaj:CB} to be a map $\psi\colon
C\o A\to A\o C$ subject to certain axioms that, again, are
sufficient to define a notion of Hopf module in $\HM.C..A.$, which
is now called an entwined module. An entwining between $C$ and $A$
naturally arises from every $C$-Galois extension $A$, in such a
way that $A$ itself is an entwined module. We will develop many of
our criteria for Galois-type extensions for entwinings between a
coalgebra $C$ and an algebra $A$ for which $A$ itself is an
entwined module, in most cases with a comodule structure induced
by a distinguished grouplike in $C$. We are most interested in the
special case where $A$ is an $H$-comodule algebra,  $C=Q$ is a
quotient coalgebra and right module of $H$, the distinguished
grouplike $e$ is the image of $1\in H$ in $Q$, and the entwining
is given by $\psi(\ol h\o a)=a\sw 0\o \ol{ha\sw 1}$. In fact, no
example of a $C$-Galois extension that would not have this form
seems to be known to date. Even though the theory of $C$-Galois
extensions is potentially more general than that of $Q$-Galois
extensions, it is perhaps more important that the formalism of
entwinings is more elegant and transparent in some situations than
the use of an auxiliary Hopf algebra, and it may serve to make
proofs more transparent and draw attention to those instances
where a Hopf algebra in the background is truly needed for more
than formal reasons.

We believe that the results of the present paper show that the
following two conditions on a $Q$-extension are of particular
interest: Equivariant projectivity, and the property that all Hopf
modules are relative injective as comodules. The status of the
former has changed significantly by our results: At its
conception, this strong, geometrically motivated condition seemed
to single out a particularly well-behaved class among Hopf Galois
extensions, and of course among the more general $Q$-Galois
extensions. Now we know that it is shared by all faithfully flat
$H$-Galois extensions when $H$ is a Hopf algebra with bijective
antipode over a field, that is, by all those Hopf Galois
extensions that are candidates for a quantum group analog of a
classical principal fiber bundle. We also prove in
\nmref{ConnectionThm} that it is shared by all $Q$-Galois
extensions with cosemisimple $Q$ over, say, the complex numbers,
another case which is of particular importance for applications to
quantum groups. On the other hand, we know very little about when
this property is fulfilled in general. The problem seems to be as
hard for quotients of Hopf algebras (i.e.\ quantum analogs of
homogeneous spaces) as it is for general $Q$-Galois extensions
(i.e.\ quantum analogs of principal fiber bundles). We have
pointed out the second interesting property, namely that all Hopf
modules are relative injective comodules, as a powerful technical
tool in criteria for $Q$-Galois situations. Again, we have
collected results that show this property to be fulfilled quite
often, but we do not know in what generality it can be proved. And
once again, the case where the $Q$-extension in consideration
comes from a quotient map $H\to Q$ seems to be as hard as the
general case.

If $H$ is finite-dimensional over a field, then both of the
conditions on a quotient coalgebra and right module $Q$ mentioned
above are equivalent to the (in general rather stronger)
condition that $H$ be $Q$-cleft. This is not hard to see using the
equivalent characterizations of the latter condition proved by
Hoffmann, Koppinen, and Masuoka \cite{Mas:FHACS,Mas:QTHA}.
Serge Skryabin kindly gave us access to his recent 
preprint \cite{Skr:PFCA}, where he proves that $H$ is in fact 
always $Q$-cleft in the finite-dimensional case.

\section*{Preliminaries and notations}
Throughout the paper $k$ denotes a commutative base ring. All maps
are at least $k$-linear, unadorned tensor product is understood to
be over $k$, algebras, coalgebras, and Hopf algebras are over $k$.
We use $\nabla\colon A\o A\to A$, $\eta\colon k\to A$ to denote
the multiplication and unit map of an algebra $A$, and
$\mu=\mu_M\colon A\o M\to M$ to denote the structure map of an
$A$-module $M$. The category of left (resp. right) $A$-modules
will be denoted $\LMod A$ (reps.\ $\RMod A$). We write
$\Delta\colon C\to C\o C;c\mapsto c\sw 1\o c\sw 2$ and
$\epsilon\colon C\to k$ for the comultiplication and counit of a
coalgebra $C$. For a right $C$-comodule $M\in\RComod C$ we write
$\delta=\delta_M\colon M\to M\o C; m\mapsto m\sw 0\o m\sw 1$ for
its comodule structure. For left comodules we use $\delta(m)=m\swm
1\o m\sw 0$. For $M\in\RComod C$ and a left comodule $N\in\LComod
C$ we denote by $M\co CN$ the cotensor product. The antipode of a
Hopf algebra $H$ is denoted by $S$.

A left $R$-module $M$ is called relative projective if it fulfills
the following equivalent conditions: Any surjective $R$-module map
$f\colon N\to M$ which splits as a $k$-module map also splits as
an $R$-module map; The module structure map $\mu\colon R\o M\to M$
splits as an $R$-module map; $\Hom_R(M,f)$ is surjective if
$f\colon N\to N'$ is a surjective $R$-module map that is
$k$-split. Note that direct summands and direct sums of relative
projective $R$-modules are relative projective. Also, if $V$ is a
$k$-module and $M$ a relative projective $R$-module, then $M\o V$
is relative projective.

Dually, a right $C$-comodule $M$ is called relative injective if
it fulfills the following equivalent conditions: Any injective
$C$-comodule map $f\colon M\to N$ which splits as a $k$-module map
also splits as a $C$-comodule map; The comodule structure map
$\delta\colon M\to M\o C$ splits as a $C$-comodule map;
$\Hom^C(f,M)$ is surjective if $f\colon N\to N'$ is an injective
$C$-comodule map which is $k$-split. Note again that direct
summands and finite direct sums of relative injective
$C$-comodules are relative injective, as is $V\o M$ whenever $V$
is a $k$-module and $M$ is a relative injective $C$-comodule.

The notion of relative projectivity (which one should call
$k$-relative projectivity, but no other versions will occur in
this paper) is a special case of the terminology of relative
homological algebra as found in \cite[Chap.IX]{Mac:H}. The same is
true for relative injectivity, provided that $C$ is $k$-flat,
which ensures that the category $\RComod C$ is abelian to begin
with.

\section{Generalities on entwining structures}

In this section we collect some general conventions and facts on
entwinings and their relation to coalgebra Galois extensions. Most
of these can be found in \cite{Brz:MACGE}; we also refer to the
survey article \cite{BrzHaj:GTENC}.
\begin{Def}\nmlabel{Definition}{entwining}
  An {\bf entwining structure} $(A,C,\psi)$ consists of an algebra $A$,
  a coalgebra $C$, and an {\bf entwining}, that is, a map $\psi\colon
  C\o A\to A\o C$ satisfying
  \begin{align*}
    \psi(C\o\nabla)&=(\nabla\o C)(A\o\psi)(\psi\o A)\colon C\o A\o
    A\to A\o C\\
    \psi(c\o 1)&=1\o c\qquad\forall c\in C\\
    (A\o\Delta)\psi&=(\psi\o C)(C\o\psi)(\Delta\o A)\colon C\o
    A\to A\o C\o C\\
    (A\o\epsilon)\psi&=\epsilon\o A\colon C\o A\to A.
  \end{align*}
\end{Def}
\begin{Def}
  Let $(A,C,\psi)$ be an entwining structure. An {\bf entwined
  module}
  $M\in\entwmod AC\psi$ is a right $A$-module and right
  $C$-comodule such that the diagram
      $$\xymatrix{M\o A\ar[r]^-{\delta\o A}\ar[d]_\mu
        &M\o C\o A\ar[r]^-{M\o\psi}&M\o A\o C\ar[d]^{\mu\o C}\\
        M\ar[rr]^-{\delta}&&M\o C}$$
      commutes.
\end{Def}
\begin{Lem}\nmlabel{Lemma}{modcomodbypsi}
  Let $(A,C,\psi)$ be an entwining structure.
  \begin{enumerate}
    \item For any $M\in\RMod A$ we have $M\o C\in\entwmod AC\psi$
    with the obvious comodule structure and the module structure
    $$\mu_{M\o C}=\left(M\o C\o A\xrightarrow{M\o\psi}M\o A\o
    C\xrightarrow{\mu_M\o C}M\o C\right).$$
    This construction defines a right adjoint functor
    $\RMod A\to\entwmod AC\psi$ to the underlying functor.

    If also $M\in\RComod C$, then $M$ is an entwined module if and
    only if the comodule structure $\delta\colon M\to M\o C$ is an
    $A$-module map.
    \item For any $M\in\RComod C$ we have $M\o A\in\entwmod AC\psi$
    with the obvious module structure and the comodule structure
    $$\delta_{M\o A}=\left(M\o A\xrightarrow{\delta_M\o A}M\o C\o
    A\xrightarrow{M\o\psi}M\o A\o C\right).$$
    This construction defines a left adjoint functor $\RComod C\to\entwmod
    AC\psi$ to the underlying functor.

    If also $M\in\RMod A$, then $M$ is an entwined module if and
    only if the module structure $\mu\colon M\o A\to M$ is a
    $C$-comodule map.
  \end{enumerate}
\end{Lem}
\begin{Rem}
  In particular, an entwining structure $(A,C,\psi)$ gives rise to
  entwined module structures on $C\o A$ as well as $A\o C$. With these
  structures, $\psi$ is a morphism of entwined modules.

  Note that the right $A$-module structure of $A\o C$ determines $\psi$ uniquely
  through the formula $\psi(c\o a)=(1\o c)a$. Dually, $\psi$ is
  determined by the $C$-comodule structure of $C\o A$.
\end{Rem}
\begin{Def}\nmlabel{Definition}{CGalois}
  Let $C$ be a coalgebra, and let $A$ be an algebra and a
  $C$-comodule. Put $B:=\rcofix AC:=\{b\in A|\forall a\in A\colon
  \delta(ba)=ba\sw 0\o a\sw 1\}$.
  Define the canonical or Galois maps
  $\beta_0\colon A\o A\to A\o C$ and $\beta\colon A\ou BA\to A\o C$ by
  $\beta_0(x\o y)=\beta(x\o y)=xy\sw 0\o y\sw 1$.
  $A$ is called a {\bf $C$-Galois
  extension} of $B$ if $\beta$ is a bijection.
\end{Def}
\begin{Lem}\nmlabel{Lemma}{canonentw}
  Let $A$ be an algebra and a $C$-comodule. If there is an
  entwining $\psi\colon C\o A\to A\o C$ for which $A\in\entwmod
  AC\psi$, then both Galois maps $\beta_0$ and $\beta$ are
  morphisms of entwined modules.

  If $A$ is a $C$-Galois extension, then there is a unique entwining
  $\psi\colon C\o A\to A\o C$ satisfying this condition.
  It is called the
  canonical entwining associated to the $C$-Galois extension $A$,
  and is given by $\psi(c\o a)=\beta(\beta\inv(1\o c)a)$.
\end{Lem}
Note that $\delta(b)=b\delta(1)$ whenever $b\in\rcofix AC$. If $A$
is $C$-Galois, we now know that $A\o C$ is a right $A$-module in
such a way that $\delta\colon A\to A\o C$ is a right $A$-module
map. Then if $\delta(b)=b\delta(1)$ for some $b\in A$, then
$\delta(ba)=\delta(b)a=b\delta(1)a=b\delta(a)$ for all $a\in A$,
hence $b\in\rcofix AC$.
\begin{Lem}\nmlabel{Lemma}{ffeq}
  Let $(A,C,\psi)$ be an entwining structure with $C$ a flat $k$-module,
  and assume that $A$ has a
  $C$-comodule structure making it an entwined module $A\in\entwmod AC\psi$
  with the regular $A$-module structure.

  The induction functor
  $$\RMod B\ni N\mapsto N\ou BA\in\entwmod AC\psi$$
  is left adjoint to the functor of coinvariants
  $$\entwmod AC\psi\ni M\mapsto \rcofix MC:=\{m\in
  M|\delta(m)=m\delta(1)\}\in\RMod B.$$
  The following are equivalent:
  \begin{enumerate}
    \item The induction functor is an equivalence.
    \item $A$ is $C$-Galois, and faithfully flat as a left
    $B$-module.
  \end{enumerate}
\end{Lem}

\begin{Lem}
  Let $A$ be a $C$-Galois extension of $B=\rcofix AC$. Assume that
  there is a grouplike element $e\in C$ with $\delta(1)=1\o e$.
  Then $\delta(a)=\psi(e\o a)$ for all $a\in A$, where $\psi$ is
  the canonical entwining. Moreover, $\rcofix MC=\{m\in
  M|\delta(m)=m\o e\}$ for all $M\in\entwmod AC\psi$.
\end{Lem}

\begin{Cor}\nmlabel{Corollary}{combyentw}
  Let $(A,C,\psi)$ be an entwining structure, and $e\in C$ a grouplike
  element. Then $A\in\entwmod AC\psi$ with the regular right $A$-module
  structure and the comodule structure
  $\delta\colon A\to A\o C$ given by $\delta(a)=\psi(e\o a)$.
\end{Cor}
In fact view $ke$ as a $C$-comodule, identify $A=ke\o A$, and
apply \nmref{modcomodbypsi}.

\begin{Lem}
  Let $(A,C,\psi)$ be an entwining structure and $e\in C$ a
  grouplike element. Endow $A$ with the $C$-comodule structure as
  in \nmref{combyentw}. Then
  $\rcofix MC=\{m\in M|\delta(m)=m\o e\}$ for every $M\in\entwmod AC\psi$.
  Assume that $A$ is $C$-Galois. Then the entwining associated to the $C$-Galois
  extension $A$ of $B$ as in \nmref{canonentw} coincides with $\psi$.
\end{Lem}
\begin{proof}
  For any $x,y\in A$ we have $\delta(xy)=\psi(e\o xy)=(\nabla\o
  C)(A\o\psi)(\psi\o A)(e\o x\o y)=(\nabla\o C)(A\o\psi)(x\sw 0\o
  x\sw 1\o y)=x\sw 0\psi(x\sw 1\o y)$.

  If $b\in\rcofix AC$, then $\delta(b)=\delta(b\cdot
  1)=b\delta(1)=b\o e$. If, on the other hand, $\delta(b)=b\o e$,
  then $\delta(ba)=b\sw 0\psi(b\sw 1\o y)=b\psi(e\o
  y)=b\delta(y)$.

  That $\psi$ coincides with the canonical entwining is a direct
  consequence of \nmref{canonentw}.
\end{proof}
\begin{Rem}\nmlabel{Remark}{leftright}
  The definition of an entwining has an obvious asymmetry (the
  coalgebra starts out on the left and ends up on the right). We
  could call an entwining as defined above a {\em right} entwining,
  and give an analogous definition of a left entwining; all the
  results collected above will then have analogous left-right
  switched versions. We will use these freely, writing
  $\tilde\psi,\tilde\beta,\tilde\beta_0,\tilde\delta$ for left entwinings,
  and the Galois maps and comodule structures associated with
  them.

  Assume now that $(A,C,\psi)$ is a bijective entwining structure,
  by which we shall mean that the map $\psi$ is bijective.
  Then the inverse $\psi\inv\colon A\o C\to C\o A$ is a left entwining.
  If $e\in C$ is a grouplike element, we thus have both a right
  $C$-comodule structure $\delta$ and a left $C$-comodule
  structure $\delta^L$ on $A$. It turns out that the left and
  right $C$-coinvariant elements of $A$ coincide: Writing $e\colon k\to A$ for the
  map that sends $1\in k$ to $e\in A$, we see that $\rcofix AC$ is
  the equalizer of $A\o e,\psi(e\o A)\colon A\to A\o C$, while
  $\lcofix AC$ is the equalizer of $e\o A=\psi\inv\psi(e\o A)$ and
  $\psi\inv(A\o e)$. We also have left versions $\beta^L\colon
  A\ou B A\to C\o A$ and $\beta_0^L\colon A\o A\to C\o A$ of the
  Galois maps, mapping $x\o y$ to $x\swm 1\o x\sw 0y$.
  Since $\psi\beta_0^L=\beta_0$ and
  $\psi\beta^L=\beta$ by the calculation $\psi\beta_0^L(x\o
  y)=\psi(\psi\inv(x\o e)y)=\psi(C\o\nabla)(\psi\inv\o A)(x\o e\o
  y)=(\nabla\o C)(A\o\psi)(x\o e\o y)=\beta_0(x\o y)$, we see that
  $A$ is left $C$-Galois if and only if it is (right) $C$-Galois.

  Similarly, given a bijective left entwining $\tilde\psi$, we get
  a right entwining $\tilde\psi\inv$, and right Galois maps
  $\tilde\beta^R,\tilde\beta_0^R$.
\end{Rem}

\begin{Rem}\nmlabel{Remark}{Hierarchie}
  \begin{enumerate}
    \item A Doi-Koppinen datum is a triple $(H,A,C)$ consisting of
    a bialgebra $H$, a right $H$-comodule algebra $A$, and a right
    $H$-module coalgebra $C$. For every Doi-Koppinen datum, we
    have an entwining structure $(A,C,\psi)$ defined by $\psi(c\o a)=a\sw
    0\o c\cdot a\sw 1$. The entwining is bijective provided that
    $H$ has a skew antipode $S^-$ (for example, $H$ is a Hopf algebra
    with bijective antipode). The inverse of $\psi$ is then given
    by $\psi\inv(a\o c)=c\cdot S^-(a\sw 1)\o a\sw 0$.
    \item\label{Qdisk} In particular, let $H$ be a bialgebra, $A$ a right
    $H$-comodule algebra, and $Q$ a quotient coalgebra and right
    $H$-module of $H$. Then we have an entwining $(A,Q,\psi)$,
    which is bijective if $H$ has a skew antipode. Note that the
    $Q$-comodule structure of $A$ is the one given in
    \nmref{combyentw} for the grouplike $e=\ol 1\in Q$. The Galois
    maps $\beta,\beta_0$ in this case are given by $\beta(x\o y)=\beta_0(x\o
    y)=xy\sw 0\o\ol{y\sw 1}$.
    \item Let $A$ be a right $H$-comodule algebra, and $C$ a {\em
    left}
    $H$-module coalgebra. Then we can view $C^\cop$ as a right
    $H^\cop$-module coalgebra, and $A$ as a {\em left}
    $H^\cop$-comodule algebra, and hence we have a {\em left} entwining
    structure $(A^\op,C,\tilde\psi)$, with $\tilde\psi(a\o c)=a\sw 1\cdot c\o a\sw 0$.
    Entwined modules in $_{A}^{C^\cop}\entwmod
    {}{}{\tilde\psi}$ are Hopf modules in $\HM.C.A..$, that is, left
    $A$-modules and right $C$-comodules $M$ satisfying
    $\delta(am)=a\sw 0m\sw 0\o a\sw 1\cdot m\sw 1$ for all $a\in
    A$ and $m\in M$. The left entwining in this situation is
    bijective if $H$ has an antipode; we have $(\tilde\psi)\inv(c\o
    a)=a\sw 0\o S(a\sw 1)\cdot c$.
    \item \label{Q'disk}A special case arises when $C=Q'$ is a quotient
    coalgebra and left $H$-module of the bialgebra $H$, and $A$ is
    an $H$-comodule algebra. Note that in this case the left Galois
    maps
    $\tilde\beta_0\colon A\o A\to C^\cop\o A$ and $\tilde\beta\colon A\ou BA\to C^\cop\o A$
    identify, respectively, with
    $\beta'_0\colon A\o A\to A\o Q'$, and
    $\beta'\colon A\ou BA\to A\o Q'$ given by $\beta_0'(x\o y)=\beta'(x\o y)=x\sw
    0y\o x\sw 1$. If $H$ is a Hopf algebra, so the left entwining
    is bijective, we can also consider the right Galois map
    $\tilde\beta_0^R\colon A\o A\to A\o Q'$, given by
    $\tilde\beta_0^R=\tilde\psi\inv\tilde\beta_0$, or $\tilde\beta^R_0(x\o y)=
    \tilde\psi\inv(\ol{x\sw 1}\o x\sw 0y)=x\sw 0y\sw 0\o \ol{S(x\sw 1y\sw 1)x\sw
    2}
    =xy\sw
    0\o\ol{S(y\sw 1)}$.
    \item\label{QQ'} Let $H$ be a Hopf algebra with bijective antipode. Then
    quotient coalgebras and right $H$-modules $Q$ of $H$
    (i.e.~coideal right ideals $I\subset H$) and quotient
    coalgebras and left $H$-modules $Q'$ of $H$ (i.e.~coideal left
    ideals $I'$ of $H$) are in bijection via $I'=S(I)$. If $Q'$
    corresponds to $Q$, then the antipode of $H$ induces a
    coalgebra anti-isomorphism $S\colon Q\to Q'$. For a right $H$-comodule
    algebra $A$, the Galois maps $\beta_0\colon A\o A\to A\o Q$ as in \eqref{Qdisk}
    and $\tilde\beta^R_0\colon A\o A\to A\o Q'$ as in \eqref{Q'disk}
    identify along $A\o S$.
  \end{enumerate}
\end{Rem}
\begin{Rem}
  Let $(A,C,\psi)$ be an entwining structure, where $A$ is a
  finite-dimensional algebra (and $k$ is a field).

  By \cite{Sch:DKHMVEM} there exists a Doi-Koppinen data, that is,
  a bialgebra $H$, a right
  $H$-module coalgebra structure on $C$, and a right $H$-comodule
  algebra structure on $A$, such that the entwining $\psi$ has the
  form given above, i.e. $\psi(c\o a)=a\sw 0\o c\cdot a\sw 1$, where
  $A\ni a\mapsto a\sw 0\o a\sw 1\in A\o H$ denotes the $H$-comodule
  structure of $A$. If we are given a grouplike $e\in C$, then
  an $H$-module coalgebra map $\pi\colon H\to C$ is given by
  $\pi(h)=e\cdot h$. The relevant Galois map for the induced right
  $C$-comodule structure on $A$ is
  $$A\o A\ni x\o y\mapsto xy\sw 0\o\pi(y\sw 1)\in A\o C.$$
  Thus, if $A$ is $C$-Galois, then $\pi$ has to be surjective,
  and we can consider $C$ as a quotient coalgebra and right
  $H$-module of $H$.

  It is not known even in the situation where $A$ is
  finite-dimensional whether $H$ can be chosen to be a Hopf
  algebra.

  There are examples of entwining structures $(A,C,\psi)$ with
  infinite-dimensional $A$ that do not come from Doi-Koppinen data
  \cite{Sch:DKHMVEM}.
  It seems to be an open question, however,
   whether there exist $C$-Galois
  extensions whose entwining cannot be induced by a Doi-Koppinen
  data.
\end{Rem}
\section{Projective Galois extensions}\nmlabel{Section}{sec:proj}
This section contains a key result of our paper, a
characterization of (relative) projective Galois-type extensions
as those for which a canonical map is a split surjective comodule
map. The result will be applied in many ways in the subsequent
sections.

The following Lemma is a combination of the adjointness in
\nmref{modcomodbypsi} (2) with the morphism $\psi$ in $\entwmod
AC\psi$, which is assumed to be an isomorphism. Its central use in
the theory of comodule algebras goes back to a paper of Doi
\cite{Doi:AWTI}.
\begin{Lem}\nmlabel{Lemma}{DoiDoi}
  Let $A$ be an algebra, $C$ a coalgebra, and $\psi\colon C\o
  A\to A\o C$ a bijective entwining. Then for each $V\in\entwmod
  AC\psi$ we have an isomorphism
  $$\Phi\colon\RComod C(C,V)\to \entwmod AC\psi(A\o C,V)$$
  given by
  $$\Phi(\gamma)=\left(A\o C\xrightarrow{\psi\inv}C\o A\xrightarrow{\gamma\o A}V\o A\xrightarrow{\mu_V}V\right).$$
  Every surjective morphism $V\to A\o C$ that splits as a
  $C$-comodule map also splits in $\entwmod AC\psi$.
\end{Lem}
\begin{proof}
  It is easy to check that $\Phi(\gamma)$ is a morphism of
  entwined modules. The inverse of $\Phi$ is given by
  $\Phi\inv(\varphi)(c)=\varphi(1\o c)$.

  If $f\colon V\to A\o C$ is a morphism in $\entwmod AC\psi$ and
  $g\colon A\o C\to V$ satisfies $gf=\id_{A\o C}$, then define
  $g_0\colon C\to V$ by $g_0(c)=g(1\o c)$, and put $\tilde
  g=\Phi(g_0)$. Then $\tilde g$ still splits $f$, since $fg_0(c)=1\o c$, and
  hence
  $$f\tilde g\psi(c\o a)=f(g_0(c)a)=f(g_0(c))a=(1\o c)a=\psi(c\o a)$$
  for $c\in C$ and $a\in A$.
\end{proof}

\begin{Thm}\nmlabel{Theorem}{projgal}
  Let $(A,C,\psi)$ be an entwining structure, and assume $A$ has a $C$-comodule structure
  making it an entwined module $A\in\entwmod AC\psi$ with the regular $A$-module structure.
  Put $B:=\rcofix AC$.
  Consider the following statements:
  \begin{enumerate}
    \item $\beta_0\colon A\o A\to A\o C$ is surjective, and splits
    as a $C$-comodule map.
    \item
      \begin{enumerate}
        \item $\beta\colon A\ou BA\to A\o C$ is bijective.
        \item $A$ is relative projective as right $B$-module.
      \end{enumerate}
  \end{enumerate}
  Then (2) implies (1).
  \\
  If $\psi$ is bijective, then (1) implies (2)(a).
  \\
  If $\psi$ is bijective, and the obvious map $A\o B\to\rcofix{(A\o
  A)}C$ is a bijection (e.g.\ $A$ is $k$-flat), then (1) implies (2).
\end{Thm}
\begin{proof}
  (2)$\implies$(1): If $A_B$ is relative projective, the
  multiplication map $A\o B\to A$ splits in $\RMod B$. Apply the
  functor $(\leer)\ou BA\colon \RMod B\to \entwmod AC\psi$ to find
  that
  $$\beta_0=\left(A\o A\cong A\o B\ou BA\xrightarrow{\mu\ou BA}A\ou BA\xrightarrow\beta A\o C\right)$$
  splits in $\entwmod AC\psi$, and in particular as a comodule map.

  (1)$\implies$(2)(a) if $\psi$ is bijective:
  By assumption, there is a $C$-colinear
  splitting of $\beta_0$, and by \nmref{DoiDoi} it follows that there is a
  splitting in $\entwmod AC\psi$.

  To prove that $\beta$ is bijective, consider more generally the
  adjunction map $\mu_V\colon \rcofix VC\ou BA\to V$ for $V\in\entwmod AC\psi$.
  We can identify $\mu_{A\o C}$ with $\beta$, and will verify that
  $\mu_{A\o C}$ is a bijection by using functoriality of $\mu$.
  Since $A\o C$ is a direct summand of $A\o A$, we only need to check that
  $\mu_{A\o A}$ is a bijection. Tensoring a $k$-free resolution of $A$ on the right with
  $A$, we can write $A\o A$ as the cokernel of a morphism between entwined modules that
  are direct sums of copies of $A$. Thus it is finally enough
  to observe that $\mu_A$ is
  bijective, which is trivial, since $\mu_A\colon B\ou BA\to A$ is the
  canonical isomorphism.

  (1)$\implies$(2)(b) if in addition $A\o B\to\rcofix{(A\o A)}C$
  is an isomorphism: We have shown that
  $A\o C$ is a direct summand of $A\o A$ in
  $\entwmod AC\psi$. Apply the functor $\rcofix{(\leer)}C$ to deduce that
  $A$ is a direct summand of $\rcofix{(A\o A)}C=A\o B$ in $\RMod B$, hence
  relative projective.
\end{proof}

\begin{Rem}\nmlabel{Remark}{leftgeneral}
  There is a left version of \nmref{projgal} for a bijective left
  entwining $\tilde\psi$,
  concerned with the condition that the left Galois map
  $\tilde\beta_0\colon A\o A\to C\o A$ splits as a left $C$-comodule
  map. This is equivalent to the condition that the right Galois
  map $\tilde\beta^R_0\colon A\o A\to A\o C$ splits as a left $C$-comodule
  map, where the left $C$-comodule structure on the source is that
  of the left tensor factor, and the one on the right is given by
  $\tilde\delta(a\o c)=\tilde\psi(a\o c\sw 1)\o c\sw 2$.
\end{Rem}

It is clear how to specialize \nmref{projgal} to the situation of
an $H$-comodule algebra $A$ and a quotient coalgebra and right
$H$-module $Q$ of a Hopf algebra $H$ with bijective antipode. By
switching sides, we also get a version for quotient coalgebras and
left modules of $H$, which we will write down explicitly to
clarify the somewhat complicated identifications.

\begin{Cor}\nmlabel{Corollary}{projQ'gal}
  Let $H$ be a Hopf algebra, and $A$ a right $H$-comodule algebra.
  Let $Q'$ be a quotient coalgebra and
  left $H$-module of $H$. Put $B:=\rcofix A{Q'}$. Consider the following statements:
  \begin{enumerate}
    \item $\beta_0'\colon A\o A\to A\o Q', x\o y\mapsto x\sw 0y\o
    \ol{x\sw 1}$ is surjective, and splits as a right $Q'$-comodule
    map.
    \item
      \begin{enumerate}
        \item $\beta'\colon A\ou BA\to A\o Q'$ is bijective.
        \item $A$ is relative projective as left $B$-module.
      \end{enumerate}
    \setcounter{weiter}{\value{enumi}}
  \end{enumerate}
  Then (2) implies (1), and (1) implies (2)(a). If the obvious map
  $A\o B\to \rcofix{(A\o A)}{Q'}$ is bijective, then (1) implies
  (2).

  If the antipode of $H$ is bijective, and $Q$ is the quotient
  coalgebra and right module of $H$ corresponding to $Q'$, then
  (1) is equivalent to
  \begin{enumerate}
    \setcounter{enumi}{\value{weiter}}
    \item $\beta_0\colon A\o A\to A\o Q,x\o y\mapsto xy\sw
    0\o\ol{y\sw 1}$ is surjective, and splits as a left
    $Q$-comodule map; here, the left $Q$-comodule structures are given by
    \begin{align*}
      A\o A\ni x\o y&\mapsto \ol{S\inv(x\sw 1)}\o x\sw 0\o y\in
      Q\o A\o A\\
      A\o Q\ni x\o q&\mapsto q\sw 1S\inv(x\sw 1)\o x\sw 0\o q\sw
      2\in Q\o A\o Q.
    \end{align*}
  \end{enumerate}
\end{Cor}
\begin{proof}
  As discussed in  \nmref{Hierarchie} \eqref{Q'disk}, we have a
  bijective left entwining $\tilde\psi$ involving $C=(Q')^\cop$.
  Applying
  the left version of \nmref{projgal} yields the stated relations
  between (1) and (2).

  As in \nmref{leftgeneral}, (1) is equivalent to the condition
  that $\tilde\beta_0^R\colon A\o A\to A\o Q'$, given by
  $\tilde\beta_0^R(x\o y)=xy\sw 0\o \ol{S(y\sw 1)}$, is surjective
  and splits as a right $Q'$-comodule map, where the comodule
  structure on the source is that of the left tensor factor, and
  that on the target is given by
  $$\delta_{A\o Q'}:A\o Q'\ni x\o q\mapsto x\sw 0\o q\sw 2\o x\sw 1q\sw 1\in A\o Q'\o Q'.$$
  If $H$ has bijective antipode, and $Q$ corresponds to $Q'$ as in
  \nmref{Hierarchie} \eqref{QQ'}, then $\tilde\beta_0^R$
  identifies with $\beta_0$ as in (3), and the right $Q'$-comodule
  structures can be identified with the left $Q$-comodule
  structures given in (3), since
  $(A\o S\inv)\delta_{A\o Q'}(x\o S(q))=x\sw 0\o q\sw 2\o S\inv(x\sw
  1S(q\sw 1))=x\sw 0\o q\sw 2\o q\sw 1S\inv(x\sw 1)$.
\end{proof}

Most of our applications of \nmref{projgal} will rely on
additional hypotheses on $C$ or the comodule structure. However,
we can draw one very general conclusion on the behavior of the
Galois condition when we pass to quotients:
\begin{Cor}\nmlabel{Corollary}{forquotients1}
  Let $(A,C,\psi)$ be an entwining structure with $A\in\entwmod
  AC\psi$ such that $A$ is a $C$-Galois extension of $B:=\rcofix
  AC$, and a relative projective right $B$-module.

  Let $(R,D,\theta)$ be a bijective entwining,
  $\pi\colon C\to D$ a surjective
  coalgebra map, and $f\colon A\to R$ a $k$-split surjective algebra
  and $D$-comodule map such that $\theta(\pi\o f)=(f\o\pi)\psi$.

  Assume that $\pi$ splits as a right $D$-comodule map.

  Then
  $R$ is a $D$-Galois extension of $S:=\rcofix RD$, and if
  the obvious map $R\o S\to \rcofix{(R\o R)}D$ is bijective, then
  $R$ is a relative projective right $S$-module.
\end{Cor}
\begin{proof}
  The commutative diagram
  $$\xymatrix{A\o A\ar[d]_{f\o f}\ar[r]^-{\beta^{(A)}_0}
  &A\o C\ar[d]^{f\o\pi}\\
  R\o R\ar[r]^-{\beta^{(R)}_0}&R\o D}$$
  shows that the canonical map $\beta^{(R)}_0$ for the $D$-extension
  $R$ is surjective.
  By assumption and \nmref{projgal} the $C$-comodule map
  $\beta^{(A)}_0$ splits. Since $f$ splits as a $k$-module map, we see that
  $\beta^{(R)}_0$ splits as a $D$-comodule map, and the claim follows
  from \nmref{projgal}.
\end{proof}
\begin{Cor}
  Let $H$ be a $k$-flat Hopf algebra with bijective antipode, and
  $Q$ a quotient coalgebra and right module of $H$. Put
  $K:=\rcofix HQ$. The following are equivalent:
  \begin{enumerate}
    \item The surjection $H\to Q$ splits as a $Q$-comodule map.
    \item $H$ is a $Q$-Galois extension of $K$ and a relative
    projective right $K$-module.
  \end{enumerate}
\end{Cor}
\begin{proof}
  (1)$\implies$(2): We apply
  \nmref{forquotients1} with $A=R=C=H$ and $D=Q$.

  (2)$\implies$(1): By \nmref{projgal} the Galois map
  $\beta_0\colon H\o H\to H\o Q$ splits as a $Q$-comodule map.
  Looking at the diagram in the proof of \nmref{forquotients1}, we see that
  $H\o\pi\colon H\o H\to H\o Q$ splits as a $Q$-comodule map by,
  say, $t\colon H\o Q\to H\o H$.
  If we define $f\colon Q\to H$ by $f(q)=(\epsilon\o
  H)f(1\o q)$, then $f$ splits $\pi$.
\end{proof}
\section{Kreimer-Takeuchi type
theorems}\nmlabel{Section}{sec:KTTT}

In this section we will discuss a generalization of the
Kreimer-Takeuchi Theorem \cite[Thm.~1.7]{KreTak:HAGEA}, which, in
turn, is a Hopf algebraic version of a result of Grothendieck on
actions of finite group schemes \cite[III, \S 2, 6.1]{DemGab:GAI}.
Let $H$ be a Hopf algebra, and $A$ an $H$-comodule algebra such
that the Galois map $\beta\colon A\o A\to A\o H$ is surjective.

The Kreimer-Takeuchi theorem says that if $H$ is finitely
generated projective, then $A$ is an $H$-Galois extension of
$B:=\rcofix AH$, and projective as left as well as right
$B$-module. A partial generalization was proved by Beattie,
D{\u{a}}sc{\u{a}}lescu, and Raianu: If $k$ is a field, and $H$ is
co-Frobenius, it follows again that $A$ is an $H$-Galois extension
of $B$ \cite[Thm.~3.2, (ii)$\Rightarrow$(i)]{BeaDasRai:GECFHA},
and at least a flat $B$-module.

In case that $k$ is a field, we will see that both results (and
the fact that $A$ is projective over $B$ also in the case studied
in \cite{BeaDasRai:GECFHA}) follow directly from \nmref{projgal};
if $k$ is not a field, projectivity of $A$ as a $B$-module
requires a little extra work. We will prove a more general result
for entwining structures, and discuss conditions under which it
applies to $Q$-Galois extensions, with $Q$ a quotient of a Hopf
algebra $H$.

\begin{Thm}\nmlabel{Theorem}{KreTak}
  Let $(A,C,\psi)$ be a bijective entwining with $A\in\entwmod
  AC\psi$.

  Assume that the Galois map $\beta_0\colon A\o A\to A\o C$ is
  surjective.

  If $C$ is $k$-flat, and
  projective as right (left) $C$-comodule, then
  $A$ is a $C$-Galois extension of $B:=\rcofix AC$ and
  projective as right (left) $B$-module.
\end{Thm}
\begin{proof}
  We only treat the version without parentheses, which
  implies the one in parentheses
  when applied to the inverse of the entwining $\psi$.

  Note that if $A$ is a projective $k$-module, then $A\o C$ is a
  projective $C$-comodule, and hence the surjection
  $\beta_0$ splits as a comodule map. The claims then follow from
  \nmref{projgal}.

  For the general case, let $M\in\entwmod AC\psi$. We have
  isomorphisms
  $$\entwmod AC\psi(A\o C,M)\cong\entwmod AC\psi(C\o A,M)\cong
  \RComod C(C,M)$$
  induced, respectively, by $\psi$ and the adjunction in
  \nmref{modcomodbypsi} (2). Since $C$ is a projective comodule by
  assumption, it follows that $A\o C$ is a projective object in
  $\entwmod AC\psi$. Thus, the surjection $\beta_0$ splits in
  $\entwmod AC\psi$, and $A$ is $C$-Galois by \nmref{projgal}.
  Now we can continue the above chain of isomorphisms by
  $$\RMod B(A,\rcofix MC)\cong\entwmod AC\psi(A\ou BA,M)\cong\entwmod AC\psi(A\o
  C,M),$$
  using, respectively, the adjunction in \nmref{ffeq}, and the
  isomorphism $\beta$.
  Now consider a surjective right $B$-module map $f\colon B^{(I)}\to A$
  for some index set $I$. Since the unit $N\to\rcofix{(N\ou BA)}C$
  of the adjunction in \nmref{ffeq}
  is a bijection for both $N=B^{(I)}$ and $N=A$, we see that
  $f$ is the coinvariant part of the surjective morphism $f\ou BA$ in
  $\entwmod AC\psi$. Since $C$ is projective as right
  $C$-comodule the two chain of isomorphisms above show that
  $\RMod B(A,f)\colon \RMod B(A,B^{(I)})\to \RMod B(A,A)$ is
  surjective, that is, $f$ splits as a $B$-module map.
\end{proof}

\begin{Rem}
  \begin{enumerate}
    \item Assume that $C$ is a projective $k$-module. If $C$ is projective
    as left $C^*$-module, then it is projective as right $C$-comodule.
    Now assume that $C$ is finitely generated projective. Then the
    converse holds, and moreover $C$ is projective as a left
    $C^*$-module if and only if $C^*$ is injective as a right
    $C^*$-module, that is, $C^*$ is a right self-injective ring.
    It is worth noting that if $k$ is a field,
    then $C^*$ is right self-injective if and only if it is left
    self-injective \cite[Thm.\ 15.1]{Lam:LMR}; in particular the hypotheses of the
    \ref{KreTak.nme} are the same for its left-right switched
    version if $C$ is finite-dimensional over a field.
    \item Assume that $k$ is a field.
    The coalgebra $C$ is called left co-Frobenius
    if there is an injective left $C^*$-module map from $C$ to
    $C^*$.
    If $C$ is left co-Frobenius, then $C$ is projective as left
    $C^*$-module
    \cite[Prop.~5]{Lin:SC}, hence projective as a right
    $C$-comodule.
    \item A Hopf algebra $H$ over a field $k$
     is left co-Frobenius as a coalgebra if and only if it
    admits a non-zero left integral $\lambda\colon H\to k$, if and
    only if $H$ is right co-Frobenius \cite[Thm.~3]{Lin:SC}.
    In this case the antipode of $H$ is bijective.
    \cite[Prop.~2]{Rad:FCHANI}.
    \item Let us say that a Hopf algebra $H$
    which is a projective $k$-module {\bf has enough right integrals}
    if the evaluation map $H^*\o H\to k$ induces a surjection $H\o I_r(H)\to k$,
    where $I_r(H)$ denotes the space of left integrals on $H$. Thus, $H$ has
    enough right integrals if and only if there are right integrals
    $\lambda_1,\dots,\lambda_k\in H^*$ and elements
    $t_1,\dots,t_k\in H$ with $\sum\lambda_i(t_i)=1$. For example,
    $H$ has enough right integrals if there is a surjective right integral
    $\lambda\colon H\to k$.

    Now $H$ is projective as a right $H^*$-module if and only if
    $H$ has enough right integrals: One checks that if
    $\lambda_i,t_i$ are as above, then
    $$\varphi\colon H\ni h\mapsto\sum ht_i\sw 2\o \lambda_i\hitby S(h\sw 1)\in H\o H^*$$
    is an $H^*$-linear splitting of the $H^*$-module structure of
    $H$, and conversely, if $\varphi\colon H\to H\o H^*$ splits
    the module structure, then $\varphi(1)\in H\o I_r(H)$ is mapped
    to $1\in k$ under evaluation.
  \end{enumerate}
\end{Rem}

Given the well-known properties of finite Hopf algebras, and the
properties of co-Frobenius Hopf algebras discussed above,
\nmref{KreTak} contains both the Kreimer-Takeuchi theorem and its
generalization in \cite{BeaDasRai:GECFHA} as special cases, when
we apply it to the entwining coming from an $H$-comodule algebra
$A$. We will be interested in the more general situation where $A$
is an $H$-comodule algebra that we view as a $Q$-extension for a
quotient coalgebra and right $H$-module $Q$ of $H$:
\begin{Cor}\nmlabel{Corollary}{KreTakApp}
  Let $H$ be a Hopf algebra with bijective antipode,
  $A$ an $H$-comodule algebra, and $Q$
  a quotient coalgebra and right $H$-module of $H$.

  Assume that the Galois map $\beta_0\colon A\o A\to A\o Q$ is
  surjective. Then it follows that $A$ is a $Q$-Galois extension
  of $B:=\rcofix AQ$ and a projective left $B$-module in each of
  the following cases:
  \begin{enumerate}
     \item $k$ is a field, and $H$ is finite-dimensional.
     \item $H$ is finitely generated projective over $k$, coflat as
    a right $Q$-comodule, and the surjection $H\to Q$ splits as a
    left $Q$-comodule map.
    \item $H$ has enough right integrals,
    is coflat as a right
    $Q$-comodule, and the surjection $H\to Q$ splits as a left
    $Q$-comodule map.
    \item $k$ is a field, $H$ is co-Frobenius, and faithfully
    coflat both as a left and a right $Q$-comodule.
    \item $H$ is $Q$-cleft and $Q$ is finitely generated
    projective.
    \item $k$ is a field, $H$ has cocommutative coradical, and $Q$
    is finite dimensional and of the form $Q=H/K^+H$ for a Hopf
    subalgebra $K\subset H$.
  \end{enumerate}
\end{Cor}
\begin{proof}
  We will verify in each case that $Q$ is a projective left $Q$-comodule
  (or a projective right $Q^*$-module, which is equivalent if
  $Q$ is finitely generated projective).
  Then the parenthesized version of
  \nmref{KreTak} can be applied to the entwining of $A$ and $Q$
  to prove the claim.

  As for (1), Skryabin \cite{Skr:PFCA} has proved that $Q^*$ is 
  Frobenius.

  Any of (2), (3), and (4) imply that $H$ is projective as left
  $H$-comodule: In the case that $H$ is finitely generated
  projective, this follows from the structure theorem for Hopf
  modules over a Hopf algebra, since $H$ can be
  considered as a Hopf module in $\HM.H^*..H^*.$ as in
  \cite{LarSwe:AOBFHA}. If $H$ has enough right integrals,
  then $H$ is projective as right $H^*$-module as we discussed
  above. (4) is a special case of (3): If $k$ is a field, and
  $H$ is left faithfully coflat over $Q$, then the
  surjection $H\to Q$ splits as a left $Q$-comodule map by
  \cite[1.1,1.3]{Sch:PHSAHA}.

  Now to prove the desired results on $Q$ under hypotheses (2), (3),
  or (4),  we may assume more generally that $H$ is a coalgebra that is
  projective as a left $H$-comodule, and $Q$ is a quotient
  coalgebra of $H$ so that $H$ is a coflat right $Q$-comodule, and
  the surjection $H\to Q$ splits as a left $Q$-comodule map.
  We have an isomorphism
    $$\LComod Q(H,V)\cong \LComod H(H,H\co QV)$$
    for any left $Q$-comodule $V$. Thus $H$ projective in
    $\LComod H$ and $H$ coflat as right $Q$-comodule implies that
    the functor $\LComod Q(H,\leer)$ is exact, thus $H$ is
    projective as left $Q$-comodule. Since $Q$ is a direct
    summand, it also is projective as left $Q$-comodule.

  We note that under hypothesis (2), with $k$ a field,
  $Q^*$ was proved to be self-injective
  by Hoffmann, Koppinen and Masuoka \cite[Thm.4.2]{Mas:QTHA}.

  Under the hypotheses in (5) Fischman, Montgomery, and Schneider
  \cite[Thm.4.8]{FisMonSch:FESHA}.
  show that $Q^*$ is Frobenius if $k$ is a field. Part of their technique
  still applies in the general
  case: We can consider $Q^*$ as a Hopf module in $\HM.Q..H.$ with the
  right $Q$-comodule structure dual to the regular left
  comodule structure of $Q$, and the right $H$-module structure defined by
  $(\theta h)(q)=\theta(qS(h))$ for $\theta\in Q^*$,
  $h\in H$, and $q\in Q$. Then $Q^*\in\HM.Q..H.$ by the
  calculation
  \begin{multline*}
    (\theta\sw 0h\sw 1)(q)\theta\sw 1h\sw 2
      =\theta\sw 0(q S(h\sw 1))\theta\sw 1 h\sw 2
      =q\sw 1 S(h\sw 2)\theta(q\sw 2S(h\sw 1))h\sw 3
      \\=q\sw 1\theta(q\sw 2S(h))
      =q\sw 1(\theta h)(q\sw 2)
      =(\theta h)\sw 0(q)(\theta h)\sw 1
  \end{multline*}
  By the results in \cite{MasDoi:GCCA} that we will review in
  \nmref{injifcleft} it follows that $Q^*$ is injective as right
  $Q$-comodule, and hence $Q$ is projective as right $Q^*$-module.

  Under the hypotheses in (6) the algebra $Q^*$ is again
  Frobenius,  by results of Fischman, Montgomery, and Schneider
  \cite[Cor.4.9]{FisMonSch:FESHA}.
\end{proof}

\section{Injectivity conditions}\nmlabel{Section}{sec:inj}

The following result and the remark following it characterize, in
particular, those $C$-Galois extensions, with distinguished
grouplike $e\in C$, that are relative injective comodules. For
Doi-Koppinen data the results are due to Doi \cite[Prop.3.2,
Prop.3.3]{Doi:UHM}, generalizing his result \cite[1.6]{Doi:AWTI}
for comodule algebras. The proofs for general entwinings are not
essentially more difficult.
\begin{Lem}\nmlabel{Lemma}{Doi}\nmlabel{Lemma}{ifinjthendirect}
  Let $(A,C,\psi)$ be a bijective entwining structure, and $e\in
  C$ grouplike. The following are equivalent:
  \begin{enumerate}
    \item $A$ is a relative injective $C$-comodule.
    \item There is a $C$-colinear map $\gamma\colon C\to A$ with
    $\gamma(e)=1$.
    \item There is a map $\varphi\colon A\o C\to A$ in $\entwmod
    AC\psi$ with $\varphi\delta_A=\id_A$.
  \end{enumerate}
  If these conditions are satisfied, then $B:=\rcofix AC$ is a
  direct summand of $A$ as a right $B$-module, the unit $N\to\rcofix{(N\ou BA)}C$
  of the adjunction in \nmref{ffeq} is a bijection for every right $B$-module $N$, and
  in particular $\rcofix{(V\ou B
  A)}C\cong V\o B$ for every $k$-module $V$.
\end{Lem}
\begin{proof}
  Clearly (3) implies (1). Assuming (1), there is at least a
  $C$-colinear map $\varphi\colon A\o C\to A$ with
  $\varphi\delta=\id_A$. Put $\gamma(c)=\varphi(1\o c)$ to prove
  (2). Assuming (2), put $\varphi:=\Phi(\gamma)$ as in
  \nmref{DoiDoi}. We have $\varphi\delta(a)=\varphi\psi(e\o
  a)=\nabla(\gamma\o A)(e\o a)=a$, proving (3).

  For a map $\varphi$ as in (3), we have in particular
  $\varphi(b\o e)=\varphi\delta(b)=b$ for $b\in B$, hence the
  coinvariant part $\rcofix{\varphi}{C}\colon A\to B$
  splits the inclusion $B\to A$.

  Also, if $\varphi\colon A\o C\to A$ is a colinear map that
  splits the comodule structure of $A$, then the pair of
  homomorphisms $\delta,A\o e\colon A\to A\o C$ is contractible in
  the sense dual to \cite[VI.6]{Mac:CWM}, since
  $\varphi\delta=\id_A$ and $\delta\varphi(A\o e)=(A\o
  e)\varphi(A\o e)$ by the calculation $\delta\varphi(a\o
  e)=\varphi(a\o e\sw 1)\o e\sw 2=\varphi(a\o e)\o e=(A\o
  e)\varphi(A\o e)(a)$. But the equalizer of a contractible pair
  is preserved under any functor, in particular under tensor
  product with a right $B$-module $N$.
\end{proof}
\begin{Rem}\nmlabel{Remark}{ifdirecttheninj}
  If $A$ is a $C$-Galois extension of $B:=\rcofix AC$,
  and $B$ is a direct summand of $A$ as right $B$-module,
  then $A$ is relative injective as a $C$-comodule.
\end{Rem}
\begin{proof}
  If $B$ is a direct summand of $A$ as right $B$-module, then
  $A\cong B\ou BA$ is a direct summand of $A\ou{B}A$ as a
  $C$-comodule. Since $A\ou{B}A\cong A\o C$ is relative injective,
  so is $A$.
\end{proof}
\begin{Rem}
  It is clear how to specialize the results to the important case
  where $A$ is an $H$-comodule algebra, and $C=Q$ is a quotient
  coalgebra and right $H$-module (here, the antipode of $H$ should
  be bijective to have a bijective entwining).
  We can also consider a quotient coalgebra and left module of $H$
  as in \nmref{projQ'gal}. Here, the relevant entwining of
  $C=(Q')^\cop$ and $A$ is
  bijective if $H$ is a Hopf algebra.
  As a result, if $A$ is injective as $Q'$-comodule, then there is
  a map $\varphi\colon A\o Q'\to A$ in $\HM.Q'.A..$ splitting the
  comodule structure, and in particular the left submodule
  $B\subset A$ is a direct summand. Conversely, if $A$ is $Q'$-Galois and
  the left submodule $B\subset A$ is a direct summand, then $A$ is
  injective as $Q'$-comodule.
\end{Rem}
\begin{Rem}
Consider a projective right $B$-module $A$. If $A$ contains $B$ as
a direct summand, then in particular $A$ is a generator. If $A$ is
a generator in $\RMod B$, then $A_B$ is faithfully flat. Now if
$B\subset A$ is a ring extension, then conversely, $A_B$
projective and faithfully flat implies that the right
$B$-submodule $B\subset A$ is a direct summand
\cite[2.11.29]{Row:RTI}. The results above and in \nmref{sec:proj}
characterize $C$-Galois extensions $B\subset A$ such that $A_B$ is
a projective generator (or has the equivalent properties we have
just discussed):
 Let $(A,C,\psi)$ be a bijective entwining structure, and $e\in
  C$ grouplike. Put $B=\rcofix AC$. Assume that $A$ is projective as $k$-module.
  Then $A$ is $C$-Galois and a right projective generator over $B$
  if and only if the canonical map $\beta_0\colon A\o A\to A\o C$
  is surjective and splits as a $C$-comodule map, and $A$ is an
  injective $C$-comodule. On the other hand, such extensions can
  also be characterized as those $C$-Galois extensions that are
  projective and faithfully flat  as right $B$-modules. Faithfully
  flat $C$-Galois extensions in turn can be characterized by the
  structure theorem for Hopf modules \nmref{ffeq}. Thus $A$ is a
  $C$-Galois extension and a projective generator as right
  $B$-module if and only if the induction functor $\LMod B\to\HM
  C..A..$ is an equivalence, and in addition $A$ is a projective
  right $B$-module.
\end{Rem}

If $k$ is a field, more can be said without assuming that $A_B$ is
projective:
\begin{Prop}\nmlabel{Proposition}{overfield}
  Let $(A,C,\psi)$ be a bijective entwining structure over a field
  $k$, and $e\in C$ grouplike such that $A$ is a $C$-Galois
  extension of $B:=\rcofix AC$ and a flat right $B$-module.
  The following are equivalent:
  \begin{enumerate}
    \item $A$ is a faithfully flat right $B$-module.
    \item The right $B$-submodule $B\subset A$ is a direct
    summand.
    \item $A$ is injective as $C$-comodule.
    \item $A$ is coflat as $C$-comodule.
    \item $A$ is faithfully coflat as $C$-comodule.
  \end{enumerate}
\end{Prop}
\begin{proof}
  Clearly (1) follows from (2).

  (1)$\implies$(5): Consider a left $C$-comodule $V$. We have a
  chain of isomorphisms
  $$A\ou B(A\co CV)\cong (A\ou BA)\co CV\overset{\beta\co{}V}\cong (A\o C)\co CV\cong
  A\o V,$$
  the first one using that $A$ is $B$-flat. By faithful flatness
  of $A_B$ it follows that $A$ is faithfully coflat, since $A$ is
  faithfully flat over $k$.

  (5)$\implies$(4) trivially, and (4)$\implies$(3) by a result of
  Takeuchi \cite[A.2.1]{Tak:FSF}.

  Finally (3)$\implies$(2) by \nmref{ifinjthendirect}.
\end{proof}

\begin{Cor}\nmlabel{Corollary}{forquotients2}
  Let $(A,C,\psi)$ be an entwining structure, and $e\in C$
  grouplike such that $A$ is a $C$-Galois extension of $B:=\rcofix
  AC$. Assume that $A$ is a relative projective right $B$-module,
  and the right $B$-submodule $B\subset A$ is a direct summand.

  Let $\pi\colon C\to D$ be a surjective
  coalgebra map with $\psi(\Ker\pi\o A)\subset \operatorname{Im}(A\o\Ker\pi)$.
  Assume that the induced map $\theta\colon D\o A\to A\o D$ is
  bijective.

  If $\pi$ splits as a right $D$-comodule map, and $C$ is relative
  injective as right $D$-comodule, then

  $A$ is a $D$-Galois extension of $S:=\rcofix AD$, a relative
  projective right $S$-module, and the right $S$-submodule
  $S\subset A$ is a direct summand.
\end{Cor}
\begin{proof}
  We already know from \nmref{forquotients1} that $A$ is a
  $D$-Galois extension of $S$ and a relative projective right
  $S$-module. In addition, since $B\subset A$ is a right module
  direct summand, we know from \nmref{ifdirecttheninj} that $A$ is
  a relative injective $C$-comodule, hence a direct summand of the
  $C$-comodule $A\o C$. Since $C$ is relative injective as right
  $D$-comodule, $A\o C$ and hence $A$ is a relative injective
  $D$-comodule. From \nmref{ifinjthendirect} we see that $S\subset
  A$ is a direct summand as right $S$-module.
\end{proof}
\begin{Rem}
  \begin{enumerate}
  \item Assume that $k$ is a field, and $\pi\colon C\to D$ is a
  coalgebra surjection such that $C$ is faithfully coflat as a
  right
  $D$-comodule. Then by \cite[A.2.1]{Tak:FSF} $C$ is injective as
  right $D$-comodule, and by \cite[1.1,1.3]{Sch:PHSAHA} the right
  $D$-comodule map $\pi$ splits.
  \item The most important application of the preceding Theorem
  occurs when $A$ is an $H$-Galois extension for a Hopf algebra
  $C=H$ with bijective antipode, and $D=Q$ is a right $H$-module
  coalgebra quotient of $H$.
  \end{enumerate}
\end{Rem}

Next, we will specialize our results to the case of $H$-comodule
algebras over a $k$-projective Hopf algebra $H$ with bijective
antipode. Here, the condition that the Galois map $\beta_0$ is
split already follows from the condition that $A$ is a relative
injective comodule. The reason is that in this case all Hopf
modules are relative injective comodules. This result is due to
Doi \cite{Doi:AWTI}. We formulate the next corollary for general
$C$-extensions with the property that all Hopf modules are
injective comodules. It shows that this is a powerful condition.
On the other hand, it is quite unclear when it is fulfilled,
although we will encounter such situations in later sections.

\begin{Cor}\nmlabel{Corollary}{ifallinj}
  Let $(A,C,\psi)$ be a bijective entwining, and $e\in C$ grouplike.
  Assume that every
  entwined module in $\entwmod AC\psi$ is relative injective as a
  $C$-comodule, and $C$ is a projective $k$-module.

  If $\beta_0\colon A\o A\to A\o C$ is surjective, then $A$ is a
  $C$-Galois extension of $B:=\rcofix AC$ and a relative projective
  right $B$-module, and the right $B$-submodule $B\subset A$ is a direct
  summand.
\end{Cor}
\begin{proof}
  The canonical map $\beta_0$ is a morphism of
  entwined modules, and a left $A$-module map. Since $C$ is
  projective over $k$, the left $A$-module $A\o C$ is projective,
  and $\beta_0$ splits as an $A$-module map. Since the kernel of
  $\beta_0$ is an entwined module and a $k$-direct summand,
  $\beta_0$ splits as a $C$-comodule map by assumption. The
  assertions now follow from \nmref{projgal} and \nmref{ifinjthendirect}.
\end{proof}
Except for projectivity of $A$ as a $B$-module, the following two
results are in \cite[Thm.3.5,Thm.I]{Sch:PHSAHA}, with a different
proof.
\begin{Thm}\nmlabel{Theorem}{injHGal}
  Let $H$ be a Hopf algebra with bijective antipode, and $A$ an
  $H$-comodule algebra. Put $B:=\rcofix AH$. Assume that $H$ is a
  projective $k$-module. The following are equivalent:
  \begin{enumerate}
    \item The canonical map $\beta\colon A\ou BA\to A\o H$ is
    surjective, and $A$ is a relative injective $H$-comodule.
    \item $A$ is an $H$-Galois extension of $B$, and the right
    $B$-submodule $B\subset A$ is a direct summand.
    \item $A$ is an $H$-Galois extension of $B$, and the left
    $B$-submodule $B\subset A$ is a direct summand.
  \end{enumerate}
  In this case $A$ is relative projective as left and right
  $B$-module.
\end{Thm}
\begin{proof}
  (1) implies (2) and projectivity of $A_B$ by \nmref{ifallinj},
  since by a result of Doi \cite[1.6]{Doi:AWTI}
  every Hopf module is a relative injective comodule.

  (2)$\implies$(1) by \nmref{ifdirecttheninj}.

  The equivalence of (1) and (3) follows by applying the one of
  (1) and (2) to the $H^\op$-comodule algebra $A^\op$.
\end{proof}

\begin{Thm}\nmlabel{Theorem}{injHGalfield}
  Let $H$ be a Hopf algebra with bijective antipode over a field
  $k$. Let $A$ be an $H$-comodule algebra, and put $B=\rcofix AH$.
  The following are equivalent:
  \begin{enumerate}
    \item The canonical map $\beta\colon A\ou BA\to A\o H$ is
    surjective, and $A$ is injective as an $H$-comodule.
    \item $A$ is an $H$-Galois extension of $B$ and a faithfully
    flat left $B$-module.
    \item $A$ is an $H$-Galois extension of $B$ and a faithfully
    flat right $B$-module.
    \item The induction functor $\RMod B\to\HM.H..A.$ is an
    equivalence.
  \end{enumerate}
  In this case $A$ is projective as a left and right $B$-module,
  and $B$ is a direct summand of $A$ as both left and right
  $B$-module.
\end{Thm}
\begin{proof}
This is a combination of the preceding Theorem with
\nmref{overfield} and \nmref{ffeq}.
\end{proof}
\section{Equivariant projectivity and injectivity}\nmlabel{Section}{strongconnectiontype}

\begin{Def}\nmlabel{Definition}{equivariantly}
  Let $R$ be an algebra, $C$ a coalgebra, and $V$ an
  $(R,C)$-bimodule, that is, a left $R$-module and right
  $C$-comodule such that $(rv)\sw 0\o (rv)\sw 1=rv\sw 0\o v\sw 1$
  for all $r\in R$ and $v\in V$.
  \begin{enumerate}
    \item $V$ is called {\bf $C$-equivariantly $R$-projective} (or just
    {\bf equivariantly projective} if no confusion is likely) if there
    is a $C$-colinear splitting of the module structure map $R\o
    V\to V$.
    \item $V$ is called {\bf $R$-equivariantly $C$-injective} (or just
    {\bf equivariantly injective}) if there is an $R$-linear splitting
    of the comodule structure map $V\to V\o C$.
  \end{enumerate}
\end{Def}

As an immediate consequence of the definition, an equivariantly
projective $(R,C)$-bimodule is a relative projective $R$-module.
In addition to the requirement that $R\o V\to V$ splits as an
$R$-module map, the definition of an equivariantly projective
bimodule requires that such a splitting can be chosen to be
equivariant with respect to the coaction of $C$. Dually, an
equivariantly injective bimodule is a relative injective
$C$-comodule.

We will show that in many interesting cases a $Q$-Galois extension
$A$ of $B$ is equivariantly projective (that is, $Q$-equivariantly
$R$-projective).

The property was studied first for Hopf Galois extensions in
\cite{DabGroHaj:SCCCPHGT}, see also \cite{HajMaj:PMDM}. It was
shown there that equivariant projectivity of an $H$-Galois
extension is equivalent to the existence of a so-called strong
connection. Connections and the strong connections introduced in
\cite{Haj:SCQPB} are algebraic analogs of differential-geometric
notions.

A very special class of extensions that have all the desirable
properties we have discussed so far is the class of cleft
extensions. The following Lemma collects properties proved by
Masuoka and Doi \cite{MasDoi:GCCA} for the case where $C=Q$ is a
quotient coalgebra and right module of a bialgebra $H$, and $A$ is
an $H$-comodule algebra; we use techniques from \cite{Sch:GHCP} in
the proof. The generalization to entwinings instead of comodule
algebras does not present additional problems.
\begin{Lem}\nmlabel{Lemma}{injifcleft}
  Let $(A,C,\psi)$ be an entwining structure, and $e\in C$ a
  grouplike element; put $B=\rcofix AC$.
  Assume that {\bf $A$ is $C$-cleft}, that is,
  there is a $C$-colinear convolution invertible map
  $j\colon C\to A$, and $C$ is a flat $k$-module.

  Then $A$ is $C$-Galois, equivariantly
  projective, and equivariantly injective. The induction
  functor $\RMod B\to\entwmod AC\psi$ is an equivalence, and
  every entwined module is injective as $C$-comodule.
\end{Lem}
\begin{proof}
  We can assume $j(e)=1$ without loss of generality; otherwise
  replace $j$ with $\tilde\jmath$ defined by
  $\tilde\jmath(c)=j\inv(e)j(c)$.

  Let $M\in\entwmod AC\psi$. Define $\pi_0\colon M\to M$ by
  $\pi_0(m)=m\sw 0j\inv(m\sw 1)$. We will first show that
  $\pi_0(M)\subset \rcofix MC$. To verify
  $$\delta_M(\pi_0(m\sw 0))\o m\sw 1=\pi_0(m\sw
  0)\o e\o m\sw 1$$ for $m\in M$ (from which the assertion follows
  by applying $\epsilon$ to the last tensor factor), we apply the
  bijective map
  $$T\colon M\o C\o C\ni m\o c\o d\mapsto (m\o c)j(d\sw 1)\o
  j(d\sw 2)\in M\o C\o C$$
  to both sides:
  \begin{multline*}
  T(\delta_M(\pi_0(m\sw 0))\o m\sw 1)
  =\delta_M(m\sw 0j\inv(m\sw 1))j(m\sw 2)\o m\sw 3
  \\=\delta_M(m\sw 0j\inv(m\sw 1)j(m\sw 2))\o
  m\sw 3=\delta_M(m\sw 0)\o m\sw 1
  \end{multline*}, and, using
  $(m\o e)a=ma\sw 0\o a\sw 1$ for $m\in M$ and $a\in A$,
  \begin{multline*}
  T(\pi_0(m\sw 0)\o e\o m\sw 1)
  =(\pi_0(m\sw 0)\o e)j(m\sw 1)\o m\sw 2
  \\=\pi_0(m\sw 0)j(m\sw 1)\sw 0\o j(m\sw 1)\sw 1\o m\sw 2
  \\=\pi_0(m\sw 0)j(m\sw 1)\o m\sw 2\o m\sw 3=m\sw 0\o m\sw 1\o m\sw 3
  \end{multline*}

  Now we can define $\pi\colon M\to \rcofix MC$ by
  $\pi(m)=\pi_0(m)$ for all $m\in M$. It is straightforward to
  check that $\rcofix MC\o C\ni m\o c\mapsto mj(c)$ is
  $C$-colinear and bijective with inverse $m\mapsto \pi_0(m\sw
  0)\o m\sw 1$. In particular every entwined module in $\entwmod
  AC\psi$ is a relative injective $C$-comodule. The Galois map
  $\beta\colon A\ou BA\to A\o C$ is bijective with inverse
  $\beta\inv(a\o c)=aj\inv(c\sw 1)\o j(c\sw 2)$. The $B$-linear
  and $C$-colinear map $A\ni a\mapsto \pi(a\sw 0)\o j(a\sw 1)\in B\o A$
  splits the left $B$-module structure of $A$, and the $B$-linear
  and $C$-colinear map $A\o C\ni a\o c\mapsto \pi(a)j(c)$ splits
  the $C$-comodule structure of $A$.
\end{proof}

We note for later use that one property noted for cleft extensions
in \nmref{injifcleft} holds more generally for equivariantly
injective extensions:
\begin{Rem}\nmlabel{Remark}{equithenall}
  Let $A$ be a $C$-Galois extension of $B=\rcofix AC$.
  If $A$ is equivariantly injective and faithfully flat as a left
  $B$-module, then every entwined module in $\entwmod AC\psi$ is
  a relative injective $C$-comodule.
\end{Rem}
\begin{proof}
  Let $M\in\entwmod AC\psi$. By assumption and \nmref{ffeq} we
  have $M\cong \rcofix MC\ou BA$, and if $\varphi\colon A\o C\to
  A$ is a left $B$-linear and right $C$-colinear splitting of the
  comodule structure of $A$, then $\rcofix MC\ou B\varphi$ splits
  the comodule structure of $\rcofix MC\ou BA$.
\end{proof}

We should stress that the condition that $A$ be $C$-cleft is much
more restrictive than the equivariant injectivity and projectivity
conditions. However, the \ref{injifcleft.nme} will also have applications to extensions $A$ that are not cleft.

More precisely, let $H$ be a Hopf algebra, and $Q$ a quotient
coalgebra and right module of $H$. It will turn out to be useful
for the study of a $Q$-extension (which need not be cleft)
to know that $H$ is equivariantly
$Q$-injective. If $H$ is finite-dimensional over a field, 
a recent important result of  Skryabin
\cite{Skr:PFCA} shows that $H$ is even $Q$-cleft for any quotient coalgebra
and right module $Q$. This had long been an open question; many equivalent characterizations of cleftness in this situation had been given by Masuoka \cite{Mas:FHACS} and Masuoka and Doi \cite{MasDoi:GCCA}, and the property
had been proved in interesting special cases by Masuoka \cite{Mas:CSFHA}.

The following general \ref{ConnectionLem.nme} links equivariant
injectivity and projectivity in a general bimodule.
\begin{Lem}\nmlabel{Lemma}{ConnectionLem}
 Let $R$ be an algebra, $C$ a coalgebra, and $V$ an
 $(R,C)$-bimodule.
 \begin{enumerate}
   \item If $V$ is a relative projective $R$-module and an
   $R$-equivariantly injective $C$-comodule, then $V$ is a
   $C$-equivariantly projective $R$-module.
   \item If $V$ is a relative injective $C$-comodule and a
   $C$-equivariantly projective $R$-module, then $V$ is an
   $R$-equivariantly injective $C$-comodule.
 \end{enumerate}
\end{Lem}
\begin{proof}
  We only show (1), the proof of (2) is dual.
  Since $V$ is relative
  $R$-projective, there is a left $R$-linear map $s_0 : V \to R \o
  V$ with $\mu_V s_0 = id_V.$ Define $s$ as the composition
  $$ V \xrightarrow{\delta_V} V \o C \xrightarrow{s_0 \o id_C} R \o V \o C \xrightarrow{id_R \o \varphi} R \o V.$$
  Then $s$ is left $R$-linear, and right $C$-colinear. Since
  $\varphi$ is left $R$-linear, $\mu_V (id_R \o \varphi) = \varphi
  (\mu_V \o id_C).$ Hence
  $$\mu_V s = \mu_V (id_R \o \varphi)(s_0 \o id_C) \delta_V = \varphi (\mu_V \o id_C)(s \o id_C) \delta_V = \varphi \delta_V = id_V.$$
\end{proof}

\begin{Rem}
  We will always apply \nmref{ConnectionLem} to the following situation:
  $A$ is an algebra and a $C$-comodule, and
  $B\subset A$ is a subalgebra such that $A$ is a $(B,C)$-bimodule
  (mostly even a $C$-Galois extension of $B$); there is a grouplike $e\in C$
  with $\delta(1)=1\o e$. From the
  proof of \nmref{ConnectionLem} we see that if $A$ is
  equivariantly injective, and there is a $B$-linear splitting
  $s_0\colon A\to B\o A$ of the left $B$-module structure of $A$
  that satisfies $s_0(1)=1\o 1$, then there is a $(B,C)$-bimodule
  splitting $s\colon A\to B\o A$ that satisfies $s(1)=1\o 1$. If
  $_BA$ is relative projective, and $_BB\subset {_BA}$ is a direct
  summand, then $s_0$ can in fact be chosen in this way.
\end{Rem}

As a consequence of \nmref{ConnectionLem} and our previous
results, equivariant projectivity always holds for $H$-Galois
extensions that are relative injective comodules, if $H$ has
bijective antipode. If $k$ is a field, this means in particular
that a faithfully flat $H$-Galois extensions for a Hopf algebra
$H$ with bijective antipode always admits a strong connection in
the sense of \cite{DabGroHaj:SCCCPHGT}.
\begin{Thm}\nmlabel{Theorem}{PropHopfGal}
  Let $H$ be a $k$-projective Hopf algebra,
  $A$ a right $H$-comodule algebra and $B= A^{\operatorname{co}H}.$
  Assume that $A$ is relative injective
  as a right $H$-comodule.

  Then $A$ is equivariantly injective. In particular, $A$ is
  equivariantly projective if and only if it is relative
  projective as a left $B$-module.

  If $B \subset
  A$ is an $H$-Galois extension and the antipode of $H$ is
  bijective, then $B \subset A$ is equivariantly projective.
\end{Thm}
\begin{proof}
  Applying \nmref{Doi} to the bijective entwining of $H^\op$ and
  $A^\op$ coming from the $H^\op$-comodule algebra $A^\op$ yields
  a right $H$-colinear and left $B$-linear map $\varphi\colon A\o
  H\to A$, so that $A$ is equivariantly injective. If $_BA$ is relative
  projective, then $A$ is equivariantly projective by \nmref{ConnectionLem} (1).
  If $A$ is an
  $H$-Galois extension, and $H$ has bijective antipode, then $A$
  is relative projective as left $B$-module by \nmref{injHGal}.
\end{proof}

For applications to quantum groups, the case of $C$-extensions
where $k$ is an algebraically closed field, and $C$ is a
cosemisimple coalgebra, is particularly important.
\begin{Rem}\label{coseparable}
\begin{enumerate}
\item A coalgebra $C$ with comultiplication $\Delta : C \to C \o
C$ is called {\bf coseparable} if there is a left and right
$C$-colinear map $\varphi : C \o C \to C$ with $\varphi \Delta =
id_C.$
 \item We will call a coalgebra $C$ {\bf right cosemisimple} if it is $k$-flat,
 and fulfills the following equivalent conditions:
 \begin{enumerate}
   \item Every right $C$-comodule is relative injective.
   \item Every right $C$-comodule is relative projective.
   \item If $M$ is a right $C$-comodule, and $N\subset M$ a subcomodule
   that is a direct summand as a $k$-module, then $N$ is a direct
   summand as a $C$-comodule.
 \end{enumerate}
 If $k$ is a field, this coincides with the usual definition
 \cite[Chap.XIV]{Swe:HA};
 note that a coalgebra over a field is right cosemisimple
 if and only if it is left cosemisimple.
 If $k$ is arbitrary and $C$ is finitely generated
 projective, then $C$ is right cosemisimple if and only if $C^*$ is a
 left semisimple algebra over $k$ in the sense of Hattori \cite{Hat:SACR}
\item A $k$-flat coseparable coalgebra is right and left
cosemisimple; see \ref{Propcosep}
 \item Let $C$ be a coalgebra over a field $k$. Then the
following are equivalent:
\begin{enumerate}
 \item $C$ is coseparable.
 \item $C$ is cosemisimple, and for any
 simple (hence finite-dimensional) subcoalgebra $D \subset C$, the
 dual algebra $D^*$ is separable.
\end{enumerate}
\item Any cosemisimple coalgebra over an algebraically closed
field is coseparable.
\end{enumerate}
\end{Rem}
According to part (3) of the preceding remark, every comodule over
a coseparable coalgebra is relative injective. This is generalized
and strengthened by the following observation:
\begin{Prop}\nmlabel{Proposition}{Propcosep}
Let $C$ be a coseparable coalgebra. Then any $(R,C)$-bimodule is
equivariantly injective.

In particular, any $(R,C)$-bimodule that is a relative projective
$R$-module, is equivariantly projective.
\end{Prop}
\begin{proof}
Let $\varphi\colon C\o C\to C$ be a left and right colinear map
satisfying $\varphi\Delta=\id_C$. Then, for any $C$-comodule $V$,
$$\varphi_V : V \o C \cong V \Box_C C \o C \xrightarrow{id_V \o \varphi} V \Box_C C \cong V$$
is a $C$-colinear retraction of the comodule structure of $V$. If
$V$ is an $(R,C)$-bimodule, then $\varphi_V$ is an $R$-module map.
\end{proof}
\begin{Thm}\nmlabel{Theorem}{ConnectionThm}
  Let $(A,C,\psi)$ be a bijective entwining structure, with $C$
  projective as $k$-module, and $e\in C$ a grouplike element. Put
  $B:=\rcofix AC$.

  Assume that the Galois map $\beta\colon A\ou BA\to A\o C$ is
  surjective.

  If $C$ is right cosemisimple, then $A$ is a
  $C$-Galois extension of $B$, projective as right
  $B$-module, and the right $B$-submodule $B$ is a direct
  summand.

  If $C$ is coseparable (for example, $k$ is an algebraically
  closed field and $C$ is cosemisimple), then $A$ is
  also projective as left $B$-module, $B$ is a direct summand as
  left $B$-module, and $A$ is equivariantly projective.
\end{Thm}
\begin{proof}
  Since every right $C$-comodule is relative injective, \nmref{ifallinj}
  implies that $A$ is a $C$-Galois extension of
  $B$, and a projective right $B$-module; \nmref{ifinjthendirect}
  implies that $B$ is a direct summand as right $B$-module.

  If $C$ is coseparable, then so is $C^\cop$, and by left-right
  symmetry, $A$ has the same properties as a left $B$-module.

  Also, it follows that $A$ is equivariantly projective by
  \nmref{Propcosep}.
\end{proof}

\section{Reduction to homogeneous
spaces}\nmlabel{Section}{sec:transfer}

Consider a Hopf algebra $H$, an $H$-comodule algebra $A$, and a
quotient coalgebra and right module $Q$ of $H$. Put $B=\rcofix AQ$
and $K=\rcofix HQ$. In this section we will collect some results
that allow us to draw conclusions on the structure of the
$Q$-extension $B\subset A$ from assumptions on the structure of
the $Q$-extension $K\subset H$. This shows that among the Galois
type extensions, which are quantum group analogs of principal
fiber bundles, the analogs of homogeneous spaces play a
distinguished role. We have already mentioned above the recent 
result of Skryabin proving that $H$ is $Q$-cleft whenever $H$ 
is finite-dimensional over a field. In particular, the hypotheses
on $H$ as a $Q$-extension of $K$ in \nmref{forHthenforA}, \nmref{Doi2},
and \nmref{QGalois} below are satisfied if $H$ is finite-dimensional, since 
these hypotheses are satisfied in the cleft case by \nmref{injifcleft}.

First, we will study the question when $A$ is equivariantly
projective. So far, we have settled this for $H$-Galois
extensions, and for the case where $Q$ is coseparable. If $k$ is a
field, $H$ has bijective antipode, $A$ is $H$-Galois and an
injective $H$-comodule, and $H$ is left and right faithfully
coflat over $Q$, we know from \nmref{forquotients2} and its
left-right switched version that
$A$ is left and right projective over the algebra $B$ of
$Q$-coinvariant elements of $A$. But we do not know whether $B
\subset A$ is equivariantly projective. As a particular case, the
following result will show that $B\subset A$ is equivariantly
projective if we assume that $H$ is $K$-equivariantly
$Q$-injective. At the same time we should stress that we do not
know when $H$ has these properties.

\begin{Thm}\nmlabel{Theorem}{forHthenforA}
  Let $H$ be a Hopf algebra, $Q$ a quotient coalgebra and right $H$-module of
  $H$, and $A$ an $H$-comodule algebra. Put $K=\rcofix HQ$, and $B:=\rcofix AQ$.

  Assume that the $(K,Q)$-bimodule $H$ is equivariantly injective, and
  that $A$ is relative injective as an $H$-comodule.

  Then the $(B,Q)$-bimodule $A$ is equivariantly injective.
  (In particular, if $A$ is a projective left
  $B$-module, then it is equivariantly projective.)
\end{Thm}
\begin{proof}
  By \nmref{Doi} there is a left $A$-linear (in particular $B$-linear)
  and right $H$-colinear map $\varphi_A\colon A\o H\to A$ with
  $\varphi_A\delta_A=\id_A$, where $\delta_A$ is the $H$-comodule structure
  of $A$. Let $\varphi_H\colon H\o Q\to H$ be a
  left $K$-linear and right $Q$-colinear map with
  $\varphi_H\delta_H=\id_H$. Now define
  $$\tilde\varphi:=\left(A\o Q\xrightarrow{\delta_A\o Q}A\o H\o
  Q\xrightarrow{A\o\varphi_H}A\o
  H\xrightarrow{\varphi_A}A\right).$$
  Then $\tilde\varphi$ is $Q$-colinear. Since $\delta_A(B)\subset B\o K$,
  $\tilde\varphi$ is also left $B$-linear. Finally
  $\tilde\varphi\ol\delta_A=\id_A$ for $\ol\delta_A$ the $Q$-comodule structure
  of $A$, since $\varphi(a\sw 0\o
  \ol{a\sw 1})=\varphi_A(a\sw 0\o \varphi_H(a\sw 1\o\ol{a\sw
  2}))=\varphi_A(a\sw 0\o a\sw 1)=a$ for all $a\in A$.
\end{proof}

Next, we will return to the criterion \nmref{ifallinj}, which says
that a surjective Galois map splits, provided every Hopf module is
a relative injective comodule. Again, we do not know in general
when this property holds, but we will see that assuming it for $H$
instead of $A$ will help.

\begin{Prop}\nmlabel{Proposition}{Doi2}
  Let $H$ be a Hopf algebra, $A$ a right
  $H$-comodule algebra that is a relative injective $H$-comodule,
  and $Q$ a quotient coalgebra and right $H$-module of $H$.

  If every Hopf module in $\HM.Q..H.$ is a relative injective
  $Q$-comodule, then every Hopf module in $\HM.Q..A.$ is a
  relative injective $Q$-comodule.
\end{Prop}
\begin{proof}
Let $M\in\HM.Q..A.$. The $Q$-colinear multiplication map
$\mu\colon M\o A\to M$ splits as a $Q$-comodule map by $M\ni
m\mapsto m\o 1\in M\o A$. By assumption the comodule structure
$\delta\colon A\to A\o H$ splits as an $H$-comodule map. Hence,
the $Q$-comodule $M$ is a direct summand of $M\o A\o H$, and it is
sufficient to check that the diagonal comodule $V\o H$ is a
relative injective $Q$-comodule for every $Q$-comodule $V$. But
$V\o H\in\HM.Q..H.$ with the $H$-module structure defined on the
right tensor factor.
\end{proof}

Note that by \nmref{equithenall}, the property required of $H$ in
\nmref{Doi2} holds in particular if $H$ is a faithfully flat
$Q$-Galois extension and equivariantly projective. This is true in
particular (or directly by \nmref{injifcleft}) if $H$ is
$Q$-cleft. In particular, the following result gives strong
conclusions on $A$ (which need not be cleft) if $H$ is $Q$-cleft.
This result is stronger than \nmref{forquotients2} in that it does
not assume that $A$ is $H$-Galois.
\begin{Thm}\nmlabel{Theorem}{QGalois}
  Let $H$ be a $k$-flat Hopf algebra with bijective antipode and $Q$ a
  quotient coalgebra and right $H$-module of $H$ such that
  $H$ is a $Q$-Galois extension of $K:=\rcofix HQ$, and a faithfully
  flat left $K$-module.

  Assume that $H$ is $K$-equivariantly $Q$-injective.

  Then every Hopf module in $\HM.Q..H.$ is a relative injective
  $Q$-comodule.

  In particular, if $A$ is a right $H$-comodule algebra which is
  a relative injective $H$-comodule, and the Galois map $A\o A\to A\o Q$
  is onto, then $A$ is a $Q$-Galois extension of $B:=\rcofix AQ$,
  a projective right $B$-module, and $B\subset A$ is a right
  $B$-direct summand.
\end{Thm}
\begin{proof}
  Every Hopf module in $\HM.Q..H.$ is a relative injective
  comodule by
  \nmref{equithenall}, by \nmref{Doi2} it follows that every Hopf
  module in $\HM.Q..A.$ is a relative injective comodule, and the
  remaining assertions follow from \nmref{ifallinj}.
\end{proof}
\begin{Rem}
  The hypothesis in \nmref{QGalois} that $H$ be a faithfully
  flat $Q$-Galois extension of $K$ is fulfilled if we assume
  (in addition to $H$ being relative injective as a
  right $Q$-comodule)
  that $H$ is left
  faithfully coflat for $\HM.Q..H.$, that is, cotensor product
  with $H$ over $Q$ preserves and reflects exact sequences in the
  category $\HM.Q..H.$.
\end{Rem}
\begin{proof}
  We vary arguments from \cite[Sec.1]{Sch:NBTCPHA}:
  We first observe that
  $\tilde\beta\colon K\o H\to H\co QH$ given by
  $\tilde\beta(x\o h)=xh\sw 1\o h\sw 2$ is an isomorphism with
  inverse $\tilde\beta\inv(g\o h)=gS(h\sw 1)\o h\sw 2$. To show
  that $\mu\colon\rcofix MQ\ou KH\to M$ is an isomorphism, it is
  enough, by hypothesis, to show that $\mu\co QH$ is an
  isomorphism. But the composition
  $$\rcofix MQ\o H\cong \rcofix MQ\ou KK\o H\xrightarrow{\id\ou
  K\tilde\beta}\rcofix MQ\ou K(H\co QH)\cong (\rcofix MQ\ou KH)\co
  QH\xrightarrow{\mu\co QH}M\co QH$$
  is given by $m\o h\mapsto mh\sw 1\o h\sw 2$; it is a morphism of
  Hopf modules in $\HM.H..H.$, and thus it is sufficient to
  observe that its coinvariant part is the identity on $\rcofix
  MQ$. That the adjunction map $N\to\rcofix{(N\ou BA)}C$ is an
  isomorphism for every $N\in\RMod B$ was already observed in
  \nmref{Doi}.
\end{proof}

\bibliographystyle{acm}
\bibliography{Promo}
\end{document}